\documentclass[aos,MSNbibl,nameyear,seceqn,dvips]{arximspdf}

\doi{10.1214/11-AOS881}
\volume{39}
\issue{5}
\pubyear{2011}
\firstpage{2266}
\lastpage{2279}

\makeatletter
\def\bsuffix #1{#1}
\makeatother

\begin{document}
\begin{frontmatter}

\title{Remembering Erich Lehmann\thanksref{T1}}
\runtitle{Remembering Erich Lehmann}

\begin{aug}
\author{\fnms{Willem R.} \snm{van Zwet}\corref{}\ead[label=e1]{vanzwet@math.leidenuniv.nl}}
\runauthor{W. R. van Zwet}
\affiliation{Leiden University}
\address{P.O. Box 9512\\
2300 RA Leiden\\
The Netherlands\\
\printead{e1}} %adresu isvedimo komanda gale!
\end{aug}
\thankstext{T1}{Supported in part by the Alexander von Humboldt
Stiftung through a Humboldt Research Award.}

% HISTORY:
\received{\smonth{12} \syear{2010}}
\revised{\smonth{1} \syear{2011}}

% ABSTRACT
%
\begin{abstract}
In this paper I shall try to sketch some typical aspects of
Erich Lehmann's contributions to statistics through his research, his
teaching, his service to the profession and his personality.
\end{abstract}

% KEYWORDS
%
\begin{keyword}[class=AMS]
\kwd[Primary ]{01A70}
\kwd{62G10}
\kwd[; secondary ]{62G20}.
\end{keyword}
\begin{keyword}
\kwd{Lehmann alternatives}
\kwd{efficiency and deficiency}
\kwd{Testing Statistical Hypotheses}
\kwd{Annals of Mathematical Statistics}
\kwd{European Meetings of Statisticians}.
\end{keyword}

\end{frontmatter}

%s1 ###
\section{Introduction}
Erich Leo Lehmann was born in Strasbourg on November~20, 1917. He was
raised in Frankfurt am Main where his father was a prominent lawyer.
Immediately after Hitler came to power in 1933, the family left Germany
for Switzerland, where Erich finished high school and studied
mathematics for two years at the University of Z\"{u}rich. In 1938 he
continued his studies at Trinity College, Cambridge and in January 1941
he arrived in Berkeley as a graduate student at the University of
California. In 1942 he moved from pure mathematics to statistics. In
1944--1945 he spent a year on Guam as an operations analyst with the
U.S. Air Force, together with his lifelong friend Joseph Hodges, Jr.
Having returned to Berkeley at the end of the war, Erich received his
Ph.D. in 1946 and after that remained on the Berkeley faculty
throughout his career. He received many honors such as membership of
both the National Academy of Sciences and the American Academy of Arts
and Sciences, as well as honorary doctorates at Leiden and Chicago. He
passed away on September 12, 2009 at the age of 91. Those who knew him
well, remember him with affection.

The purpose of this paper is not to present an authoritative account
of Erich Lehmann's many achievements. That will doubtless be done at
some later date. My aim is merely to describe, from a purely personal
point of view, the qualities that impressed me most about Erich Lehmann
as a scientist and a person, and that influenced my own life to a
considerable extent.

Luckily there exists an excellent book about Erich's life and times
written by himself [\citet{Leh08}]. Although this book is formally
about his interactions with fellow scientists, it does present a
reasonably complete picture of his thought and activities over the
years. An earlier version of some of the material in this book may be
found in an interview that appeared in \textit{Statistical Science}
[\citet{DeG86}]. In \citet{Leh97} Erich provided an account
of the
history of his famous book \textit{Testing Statistical Hypotheses}.
Another rich source of material is Constance Reid's biography of Jerzy
Neyman [\citet{Rei82}]. Finally there are my own conversations and
correspondence with Erich that intermittently cover a period of roughly
45 years. I shall use these sources freely to discuss some of Erich's
activities and thoughts.

%s2 ###
\section{Research}
Erich's research was paramount to his person, so I shall begin by
discussing two of his papers that I particularly like and that are
typical of his research style. Both papers start out with deceptively
simple ideas that yield results that greatly improve our understanding
of essentially complex matters. A more comprehensive view of Erich's
research output will be available in a volume of his selected papers
that should appear shortly.

%s2.1 ###
\subsection{The power of rank tests}
Lehmann
(\citeyear{Leh53}) is the first of
these two papers. I discussed it for the volume of Lehmann's selected
papers, but I think the paper has a natural place here, too. It
consists of two simple remarks, but of course these remarks were only
found to be simple after they had been made! The setting is the power
computation for two-sample rank tests. Let $X_1,\ldots ,X_m, Y_1,\ldots
,Y_n$ be independent, with the $X_i$ and the $Y_j$ having continuous
distribution functions $F$ and $G$, respectively. We wish to test the
hypothesis $F=G$ against the alternative that the $Y$'s are
stochastically larger than the $X$'s, that is, $F(x)\geq G(x)$ for all
$x$ with strict inequality for at least one value of $x$ and hence on
an interval. We rank the combined sample in increasing order, let
$R_1,\ldots ,R_n$ denote the ranks of $Y_1,\ldots ,Y_n$ and define the rank
test statistic $T=\sum_{1\leq j \leq n} k(R_j)$ where $k$ is an
increasing function. The rank test rejects the hypothesis for large
values of $T$. The test is distribution-free in the sense that under
the hypothesis the distribution of $T$ is independent of the common
distribution $F$ of the observations. Of course this is immediately
obvious since under the hypothesis the elements of the combined sample
of $X$'s and $Y$'s are independent and identically distributed (i.i.d.)
and therefore exchangeable, and hence every permutation of the ranks of
the combined sample is equally probable. There is also a slightly
different and perhaps less obvious argument for showing that the test
is distribution-free. One simply notes that by subjecting all
observations to the same increasing transformation, the common
distribution $F$ of the combined data may be changed into any other
distribution, but the ranks---and hence their distribution---remain
unchanged. As we shall see, the advantage of the second argument is
that it does not involve the exchangeability of the combined sample
under the hypothesis and the explicit knowledge of the joint
distribution of the ranks that follows from this.

The problem with rank tests was the computation of the power, that is,
the probability of rejection under an alternative. Concentrating on
shift alternatives $G(x)=F(x-\theta)$ for $\theta>0$, it took the
combined efforts of the statistics profession quite some time before
this question found a satisfactory asymptotic answer as $m$ and $n \to
\infty$. For fixed sample sizes the problem was considered hopeless in
general. Lehmann noticed that under the alternative, part of the
invariance that had proved so useful under the hypothesis still
survives. An increasing transformation of all of the observations still
leaves the ranks unchanged. So if $K$ is a distribution function on
$(0,1)$ and one considers the alternative $(F,K(F))$, that is,
$G(x)=K(F(x))$, then the increasing transformation $F$ of the random
variables carries this into the alternative $(U,K)$ where $U$ denotes
the uniform distribution on $(0,1)$ for the $F(X_i)$ and $K$ the
distribution function of the $F(Y_j)$. Hence for every $F$, the
distribution of $T$ under the alternative $(F,K(F))$ is the same as
that under $(U,K)$.

The second, even simpler remark was that two years earlier Wassily
Hoeffding had derived an expression for the distribution of the ranks
in terms of the distribution of uniform order statistics [Hoeffding
(\citeyear{H1951})]. This expression is particularly easy to compute explicitly for
the case where $K(v)=v^a$ for $a>1$. So Lehmann suggested considering
the composite alternative $(F,F^a)$ for fixed $a>1$ and all continuous
$F$, for which the distribution of the ranks---and hence of $T$---is
independent of $F$ and equal to that for $(U,U^a)$ which is easy to
compute for any sample size. As a result one can now compute the exact
small sample power of the rank test for a curve of alternatives
$(F,F^a)$ for varying $F$ and limiting results for $m, n \to\infty$
are readily available. With the aid of Hoeffding's formula for the
distribution of the ranks under the alternative and the Neyman--Pearson
lemma, one can find most powerful and locally most powerful tests for
such alternatives.

Of course this does not solve the problem of computing the power
against shift alternatives because in that case the computation is
still unpleasant and the resulting power will depend on $F$ as well as
on the shift $\theta$. In defense of his $(F,K(F))$ models for varying
$F$ as opposed to the usual shift models, Lehmann points out that it
seems that when distribution-free methods are appropriate ``\textit{one
usually does not have very precise knowledge of the alternative. What
is then required are alternatives representative of the principal types
of deviation from the hypothesis, in terms of which one can study, at
least in outline, the ability of various tests to detect such
deviations}.'' Some 60 years later many of us would agree with this sentiment.

During a meeting in Berkeley on the occasion of Erich's 80th birthday
in 1997, I spoke about this paper and evoked the image of a large group
of statisticians wearing themselves out by climbing the mountain of
technicalities associated with the study of shift alternatives, while
at the same time Erich Lehmann tiptoed quietly around the mountain to
greet them when they got to the other side. What I did not realize at
the time and only discovered when re-reading Erich's interview in \textit{Statistical Science} [\citet{DeG86}] was that Erich himself also liked
this paper particularly and considered it one of his favorites. After
the meeting he wrote to me: ``\textit{It is rather a strange experience---at the same time embarrassing and ego-swelling---to find one's work the
topic of several lectures. I hope it does not sound conceited if I say
that your talk was the one that gave my work some `personality'}.'' It
seems safe to think that he liked the meeting as well as the interest
in his 44-year-old paper.

%s2.2 ###
\subsection{Deficiency}
Hodges and Lehmann (\citeyear{HodLeh70}) is the other paper I
would like to discuss. Again the paper advances a simple idea that
could almost be called obvious, except for the fact that nobody had
considered it in any generality before. Earlier studies of this kind,
such as \citet{Fis25} or \citet{Rao61}, are concerned with a single
specific problem in a restricted setting. The general approach set
forth in the Hodges--Lehmann paper led to a revival of asymptotic
expansions of distribution functions that substantially strengthened
large sample theory. Also the paper has the unique distinction of
having a title consisting of a single word!

To compare the quality of two statistical procedures $A$ and $B$ for
the same problem, one may compare the sample sizes $N_A$ and $N_B$
needed to obtain the same performance from the two procedures. As a
measure of performance one can think of the variance $\sigma^2$ of an
estimator or the power $\pi$ of a test of a null hypothesis against a
given alternative for a given level of significance $\alpha$. If we
restrict attention to tests rather than estimators, $N_A$ and $N_B$ are
the sample sizes needed to obtain a power $\pi\in(0,1)$ against the
given alternative for a significance level $\alpha\in(0,1)$ with
procedure $A$ or $B$. The ratio $N_B/N_A$ may then be viewed as the
efficiency of procedure $A$ with respect to $B$, that is, as a measure
of how well procedure $A$ performs for this particular testing problem
when compared to $B$. In most cases of interest, however, the quantity
$N_B/N_A$ is too difficult to compute, so one resorts to asymptotics.

Let $X_1, X_2,\ldots , X_n$ be i.i.d. random variables with common
density $f_\theta$ and consider the problem of testing the null
hypothesis $H_0 \dvtx  \theta=0$ against the alternative $K(c,n) \dvtx  \theta
=cn^{-1/2}$ for $c>0$. For sample sizes $n=1,2,\ldots $ we thus have a
sequence of testing problems and under the usual regularity conditions
the alternative hypothesis approaches the null hypothesis as $n \to
\infty$ in the sense of mutual contiguity. The level of significance
$\alpha\in(0,1)$ is kept fixed and contiguity ensures that the power
of any sequence of tests against this sequence of alternatives is
bounded away from 0 and 1 as $n \to\infty$. Suppose, as is often the
case, that as $n \to\infty$ the power for the $n$th testing problem
against $K(c,n)$ equals
%
%e2.1 ###
\begin{eqnarray}\label{1.1}
\pi_{n,A}(K(c,n)) &=& a_0(c) + a_1(c)n^{-1/2} + a_2(c)n^{-1}+ o(n^{-1}),
\nonumber
\\[-8pt]
\\[-8pt]
\pi_{n,B}(K(c,n)) &=& b_0(c) + b_1(c)n^{-1/2} + b_2(c)n^{-1}+ o(n^{-1})
\nonumber
\end{eqnarray}
for test procedures $A$ and $B$, respectively, smooth functions $a_i$
and $b_i$ and $c>0$. If procedure $A$ is applied for sample size $n$,
then procedure $B$ should be applied with sample size $k=k_n$ in order
to obtain the same power, so
%
%e2.2 ###
\begin{equation}\label{1.2}
\pi_{n,A}(K(c,n)) = \pi_{k,B}(K(c,n)) = \pi_{k,B}\bigl(K\bigl(c (k/n)^{1/2},k\bigr)\bigr)
\end{equation}
and hence
%
%e2.3 ###
\begin{equation}
a_0(c) = b_0\bigl(c (k/n)^{1/2}\bigr) + o(1).
\end{equation}
Typically, the functions $a_i$ and $b_i$ are smooth and this will allow
us to find the \textit{asymptotic relative efficiency} of procedure $A$
with respect to $B$ for the sequence of alternatives $K(c,n)$
\[
E(A,B) = \lim_{n \to\infty} \frac{N_B}{N_A} = \lim_{n \to\infty} \frac{k_n}{n}
\]
as the solution $E$ of the equation
%
%e2.4 ###
\begin{equation}
a_0(c) = b_0(cE^{1/2}).
\end{equation}
If $E(A,B) >1$, then this indicates that for large samples procedure
$A$ is preferable to procedure $B$ for testing against the alternative
$K(c,n)$ as procedure $B$ needs a sample of size $nE(A,B)$ whereas
sample size $n$ suffices for procedure $A$ to achieve the same power.

The problem becomes more interesting if the limiting powers $a_0(c)$
and $b_0(c)$ are equal and hence $E(A,B)=1$. In this case the
asymptotic relative efficiency $E$ provides no guidance for the choice
of one of the procedures, but so far we have only used the leading
terms $a_0(c)$ and $b_0(c)$ in (\ref{1.1}). Hodges and Lehmann
suggested substituting the full expansions (\ref{1.1}) in (\ref{1.2})
and using the fact that $k/n \to1$ to expand
\[
\biggl (\frac{k}{n} \biggr)^{1/2} =  \biggl(1+\frac{k-n}{n} \biggr)^{1/2}= 1
+ \frac{k-n}{2n} - \frac{(k-n)^2}{8n^2} + \cdots
\]
and use the smoothness of the functions $b_i$ to solve (\ref{1.2}) and
obtain an expansion for $(k-n)/n$. They call $d_n = (k_n-n)$ the
deficiency of procedure $B$ with respect to procedure $A$ and when
$E(A,B)=1$ one typically obtains an expansion of the form $d_n = h_1
n^{1/2} + h_2 + o(1)$. The case where $h_1 =0$ so that procedure $B$
needs only a bounded number of additional observations to compete with
procedure $A$ is of course of particular interest. This occurs if
$a_1(c) = b_1(c)$. Note that our informal discussion does not cover
some of the more complex cases where the functions $a_i(c)$ and/or
$b_i(c)$ may not be bounded and deficiencies of different orders of
magnitude than $n^{1/2}$ or 1 may occur.

In their paper Hodges and Lehmann computed deficiencies for a number
of pairs of parametric tests and estimators. The paper ends with a
section ominously entitled ``Further possibilities.'' The authors write:
``\textit{More interesting perhaps are a number of problems in which the
deficiency concept appears to be useful, but where its application
presents certain technical difficulties stemming from the fact that the
computation of deficiency requires higher-order asymptotic terms than
we encounter in the usual efficiency analysis}.'' Of course this refers
to the fact that power expansions like (\ref{1.1}) are needed for
deficiency calculations, whereas the leading terms $a_0(c)$ and
$b_0(c)$ suffice for computing efficiencies. They continue by
describing a few of these problems of which the first one is simply:
``\textit{What is the deficiency (for contiguous normal shift
alternatives) of the normal scores test or of van der Waerden's $X$-test
with respect to the $t$-test?}'' After their earlier work on rank tests
they clearly expected the rank tests to do well.

When reading the paper, both Peter Bickel in Berkeley and I in Leiden
understood that here was our homework assignment for the next few
years. In the discussion of \citet{Leh53} above, I mentioned the
mountain of technicalities involved in computing the limiting power of
rank tests for location alternatives, and now we were asked to find
asymptotic expansions for these powers with remainder $o(n^{-1})$. All
I knew about this subject was the chapter in \citet{Fel66} on Edgeworth
expansions for sums of i.i.d. random variables and that did not seem
quite sufficient. However, after months of hard labor it turned out
that Hodges and Lehmann's optimistic view of rank tests was right. We
found that for contiguous normal location alternatives, the deficiency
of the normal scores test or van der Waerden's test with respect to the
$t$-test equals $d_n = (1/2)\log\log n + O(1)$ for the one-sample case
and $d_n = \log\log n + O(1)$ for the two-sample case. For large $n$,
only $\log\log n$ additional observations are needed to make the best
rank tests perform as well as the best parametric test! Also, the
normal distribution is the slightly unpleasant exception referred to
above. A more typical result is that for contiguous logistic location
alternatives, the deficiency of Wilcoxon's test with respect to the
most powerful parametric test tends to a finite limit for the one- as
well as the two-sample problem. Even better news was that for
contiguous normal location alternatives the deficiency of the
permutation test based on sample means with respect to the $t$-test tends
to zero as $n \to\infty$. Explicit formulas for these deficiencies may
be found in \citet{AlbBicvan76} and \citet{Bicvan78}.

When Frank Wilcoxon introduced his one- and two-sample rank tests in
\citet{Wil45}, the general feeling---and to some extent his own---was that these rank tests were quick-and-dirty methods that would
suffer from a serious loss of power compared to classical parametric
tests. Since then research has consistently shown that this is not the
case. Through his own work and by pointing out the right problems to
others, Erich Lehmann played a major role in this development.

%s3 ###
\section{Teacher and friend of many}
The Latin text on the diploma of Erich Lehmann's honorary doctorate at
Leiden praises him for contributing substantially to the formulation of
statistical science as a consistent mathematical theory, and for his
work in areas such as nonparametric methods, estimation theory, robust
methods and second-order asymptotics. The text then continues
describing Erich also as ``\textit{magister et amicus multorum},'' that is,
teacher and friend of many, and adds that through the books he has
written, he has taught and educated scholars all over the world. Though
all of this was written in somewhat stilted Latin, it does give a
succinct description of Erich Lehmann's accomplishments. Two typical
examples of his research in nonparametrics and second-order asymptotics
were discussed in the previous section. Let us now discuss Erich's book
\textit{Testing Statistical Hypotheses} and its role in the
mathematization of statistics and the education of an entire generation
of statisticians worldwide.

In \citet{Leh08} Erich tells how he came to Berkeley in January 1941
to study pure mathematics, and after his first semester got a teaching
assistantship with Evans, chair of the mathematics department. After
America's entry into the war, the university ran a training program for
the army with Erich as one of the instructors. In the summer of 1942
Evans advised him that it might be more useful for the war effort if he
would switch from pure mathematics to statistics with Jerzy Neyman. In
the fall of 1942 Erich took Neyman's first upper division course in
statistics. At some time during the course Neyman suddenly had to leave
for three weeks and Erich had to take over the lectures that were
supposed to be his introduction to statistics as a student. To make
things even scarier, there was no text for the course.

After finishing the first semester of the statistics program, Erich
decided that he did not like statistics. He and other pure mathematics
students thought statistics did not possess the beauty of number theory
or other parts of pure mathematics. He felt that ``\textit{ad hoc methods
based on questionable and quite arbitrary assumptions were used to
solve messy problems}.'' He went to the newly appointed famous logician
Alfred Tarski and asked if he could study algebra with him. Like
Neyman, Tarski was Polish and intent to build up his own group in
Berkeley. That was where the similarity of the two ended and in
Berkeley they were called ``\textit{Poles apart}.'' Of course Tarski was
happy to have this new student, but before Erich had a chance to tell
Neyman he was leaving, Neyman offered him a promotion to lecturer at a
much higher salary than his assistantship. Erich decided he could not
afford to turn down this offer. Had the outcome been different, then
this article would have been written by a different person for a
different journal, with the word ``statistics'' replaced by ``algebra'' throughout.

Having decided to stay in the statistics program, Erich started by
attending Neyman's basic graduate course in statistical theory and
discovered that statistics was not so bad after all. In the 1930s Jerzy
Neyman and Egon Pearson had succeeded in formulating a consistent
mathematical model that made it possible to discuss the properties of
statistical tests in mathematical terms. The course was largely based
on this work and Erich quickly saw that there was sensible mathematics
going on. In 1944--1945 there was a break in his graduate education
when he served for a year as an operations analyst with the U.S. Air
Force on Guam to study bombing accuracy. In 1946 he completed his Ph.D.
with Neyman as advisor, and Hsu and P\'{o}lya standing by during Neyman's
long absence. Many of us had less impressive advisors!

After obtaining his Ph.D., Erich's teaching activities started in
earnest. In 1948 and 1949 Neyman entrusted him with the basic graduate
course on hypothesis testing. Erich's first graduate student Colin
Blyth attended the course and took notes, and after careful reading by
Erich the notes were mimeographed and sold at cost, first by the
Statistical Laboratory and later by the University Bookstore. As it was
the only systematic treatment of the Neyman--Pearson theory, the notes
began to be used at other universities in the United States and abroad.
Of course the notes covered only the basics of the theory and the
simplest applications, so there was soon an increasing demand to expand
the material into a textbook.

At this point a serious conflict occurred. As mentioned in \citet{Rei82}
and \citet{DeG86}, Neyman heard rumors that Lehmann had not followed
the script for the course closely enough in his class. Indeed, Erich
had added some newer material, for example, on invariance. As a result,
Neyman would not let Erich teach the course anymore, and it was only
after Neyman stepped down as chair of the then new Statistics
Department in 1956, that Erich taught it again.

Of course this did not stop Erich from working on the book and ten
years after Blyth's notes the first volume \textit{Testing Statistical
Hypotheses} appeared in 1959 [\citet{Leh59}]. I first saw the book when
I started work on my thesis in Amsterdam in January 1961, and it made
an overwhelming impression on me. Of course I had been taught the basic
ideas of hypothesis testing and the derivation of the distributions of
various common test statistics, but it all seemed to consist of a lot
of special cases without much in common. Already when reading the first
chapter on decision theory things started to fall into place. As a
master student in mathematics I had taken a full year of measure theory
with Zaanen in Leiden and another year of measure theoretic probability
with van Dantzig in Amsterdam, so I could skip the probability
background in Chapter 2. The remaining chapters provided a lucid
account of the subject along the lines of Neyman and Pearson, in a
logical order and written with the utmost care.

A few years earlier I had taken a course of van Dantzig where he
discussed R.~A. Fisher's attempt at providing a framework for
statistical inference in his book \textit{Statistical Methods and
Scientific Inference} [\citet{Fis56}]. I am afraid I was unable to
understand the content of this book in any mathematical sense, and neither
was van Dantzig who wrote a scathing review of the book [\citet
{autokey3}]. With this earlier experience in mind I was relieved to find
Professor Lehmann's account of the Neyman--Pearson approach
intelligible and entirely satisfactory, even though Professor Fisher
has warned us that ``\textit{it is to be feared, therefore, that the
principles of Neyman and Pearson's `Theory of Testing Hypotheses' are
liable to mislead those who follow them into much wasted effort and
disappointment, and that its authors are not inclined to warn students
of these dangers}'' [\citet{Fis56}, page~88].

I am afraid I am one of those people whom Professor Fisher feared
would be misled, but I will admit it does not bother me. I was
certainly not the only one for whom \textit{Testing Statistical
Hypotheses} was a revelation. The claim I made above that the book
educated an entire generation of statisticians worldwide cannot be far
off the mark. Unfortunately the sales number of the book is unknown,
but in \citet{Leh97} Erich writes that the money he made was not
nearly enough to imitate Landau and build a mansion in G\"{o}ttingen
from the proceeds of the book, but would certainly have been sufficient
to buy a fancier car than the one Erich was actually driving.

Unfortunately, total perfection does not exist and in \citet{Leh97}
Erich writes: ``\textit{The biggest source of errors was the more than 200
problems. The difficulty often resulted from some fine points or
special cases that I had overlooked and that naturally caused readers
much trouble when they struggled with them. Letters asking for
clarification were not only painful reminders of my ineptitude but they
could also take quite a bit of time and effort to answer at a time when
I was no longer working in this area. I was saved from this bondage to
my past errors when a heroic group of 15 Dutch statisticians decided to
work though the whole collection systematically and in 1984 with
Wiley's permission published the solutions as a 310-page book
[{\setattribute{cite}{font}{\rmfamily}\citet{Kal84}}].  One member of the group told me later that
this was the most painful job he had ever undertaken. However, from
then on I was able to answer queries about the problems with a simple
reference}.'' I have to admit that I was not one of the 15 Dutch heroes.

The companion volume \textit{Theory of Estimation} appeared in 1983. It
also had its roots in Colin Blyth's mimeographed lecture notes dating
back to 1950, but the many developments since that time made even more
rewriting and addition necessary. In Erich's own words in the
introduction, the two volumes together ``\textit{provide an introduction
to classical statistics from a unified point of view}.'' What was still
new in the 1950s was classical in the eighties! Since then both volumes
have been revised and reprinted and there is now a third edition of
\textit{Testing Statistical Hypotheses} with Joseph Romano as a co-author
and double the number of pages, and a second edition of \textit{Theory of
Point Estimation}, joint with George Casella. Both books continue to
acquire new friends. In between Erich also published \textit{Elements of
Large Sample Theory}, a more elementary---and easier to
read---treatment of asymptotic theory than most comparable texts. His
last book on Fisher, Neyman and the creation of classical statistics
will hopefully appear in 2011. When working on this book he wrote to me
with his typical self-deprecating humor: ``\textit{I lead the kind of life
appropriate to my age: Taking pills, going to doctors, and sleeping
during the day as well as at night. So as not to be completely idle, I
putter away at a new book: Fisher, Neyman, and the creation of
classical statistics. It is a story that I find not only scientifically
interesting but also full of drama}.'' If I am allowed a single
commercial, let me say that this will be the book that all of us should read.

One last word about Erich as ``teacher and friend of many.'' Erich was
always friendly to young people. A year after getting my Ph.D., I
visited the United States for the first time and had a chance to attend
the Fifth Berkeley Symposium. When I entered the coffee room at the
symposium, a gentleman sought me out, introduced himself as Erich
Lehmann and started talking about my thesis. I could not believe my
ears: the famous Professor Lehmann taking time to speak to me and
actually having an idea what my thesis was about! This was typically
Erich: kind, considerate and generous with the younger generation,
including of course his own students. This may explain that even in his
nineties when his own generation had largely vanished, he had so many
friends and admirers in the profession. They got together at the
Lehmann symposia, spoke about recent developments and generally enjoyed
each other's company. This is probably the best recipe for staying
young, which Erich succeeded in doing until a remarkable age. I am told
that in the end he was  ready to go.
We have known a great man that we will sorely miss.

%s4 ###
\section{Administrative and governing tasks}
First of all, I should perhaps dispel a common misconception about
Erich Lehmann that he actively encouraged himself. I am referring to
his self-professed disinterest and inability in matters of an
administrative or governing nature.

Many professors exhibit this behavior in hopes of being spared
administrative jobs through assumed incompetence. However, even though
he may not have liked it, Erich knew very well how to run things and
convince people. He was chair of the Statistics Department in Berkeley
and the faculty seems to agree that he was very effective: when he
wanted something done, he was usually right and got it done. In his
obituary for Erich in \textit{Bernoulli News}, Peter Bickel phrased this more
elegantly by referring to Erich's ``\textit{great astuteness about the
world, what could be achieved, and how to do it}.''

From 1953 through 1955 Erich was editor of \textit{The Annals of
Mathematical Statistics}. In those days, the Annals did not have a
separate managing editor, so the editor's work did not stop with the
final acceptance of papers to be published. He was also responsible for
seeing the issues through the printing process. No wonder Erich would
have much preferred to continue as associate editor. However, Neyman
strongly encouraged him to accept the editorship, and was willing to
reduce his teaching load and provide salary and space for an editorial
assistant. He obviously felt that Erich would do a good job, but he
also thought that this would give his Statistical Laboratory---soon to
turn into the Department of Statistics---some additional visibility.

So Erich accepted the editorship and apparently did very well, since
after his three-year term he was asked to continue for another term.
Neyman agreed to continue the previous arrangement, Erich accepted a
second term, but some time later Neyman turned around and said that he
could no longer provide an editorial assistant. Erich found this
incredible and told IMS president Henry Scheff\'{e} about this, who
suggested that IMS could probably pay the assistant's salary. When
Erich told Neyman about this suggestion, Neyman replied that he would
also need the room that the \textit{Annals} was using for its editorial office
and would only have a windowless storeroom in the attic available for
the \textit{Annals}. Neyman argued that he needed the money and the office space
for the editorial assistant of the proceedings of the upcoming Third
Berkeley Syposium. Erich realized that further discussion was useless
and resigned as editor of the \textit{Annals}. An acrimonious exchange of
letters started between Neyman and various IMS council members, and in
the end even David Blackwell, the most reasonable person around, found
Neyman unwilling to listen. This is not the place to discuss the
motives behind this quarrel, but only to point out that Erich felt that
this amounted to a definite break between\vadjust{\goodbreak} Neyman and himself that
lasted for many years. He makes it clear, however, that Neyman never
took this personal problem out on his career but always tried to
further it. The entire affair is recalled in detail in \citet
{Rei82}. A
quarter of a century later Erich writes that he thought ``\textit{the
issue didn't warrant such an uproar and did what he could to calm the
waters}'' [\citet{Leh08}]. All in all, it was a rather sad ending
of an
editorial job well done.

Let me relate a less-known instance of Erich's effectiveness that
ended on a much happier note. Around 1960 the only statistics meetings
in Europe were the biennial sessions of the International Statistical
Institute (ISI) that were held all over the world, and with some
regularity in Europe. However, these meetings were largely devoted to
official statistics and a number of people felt that there should also
be meetings of mathematical statisticians in Europe. Statistics was
growing rapidly in Europe, the United States was far away, and without
regular international contacts people felt isolated. Jim Durbin felt
strongly about this and so did Henri Theil, president of the
Econometric Society, who wanted joint meetings in Europe of
statisticians and econometricians. Theil asked ISI to sponsor such
European meetings, but ISI declined. Durbin approached IMS. At the time
Erich was president-elect of IMS and an informal meeting to discuss
this matter was held in Erich's home in Berkeley. Erich felt strongly
that IMS as the largest scientific society in mathematical statistics
should be willing to take on an international role. He argued that one
quarter of IMS members lived outside the U.S. and that this proportion
was rising. At its meeting in Stanford in August 1960 the IMS Council
decided to initiate the holding of European regional meetings in
addition to the annual meeting and the meetings of the American
regions. A European regional committee was appointed and at the IMS
council meeting in Seattle in June 1961 with Erich as president, this
committee recommended to hold such meetings in September 1962 in Dublin
and in July 1963 in Copenhagen. Council approved.

Thus the first European meeting at Dublin took place two years after
the idea first came to IMS, which is a remarkable performance for which
Erich's support was indispensable. He recognized a great opportunity
when he saw it and acted accordingly to convince his less
internationally minded colleagues. The Dublin meeting was truly
memorable, with an excellent scientific program attended by 300 people,
and a superb party hosted by the Guinness brewery in memory of their
former brewmaster W.~S. Gosset, also known as Student. The European
meetings of statisticians continue to the present day and have played a
major role in the development of statistics in Europe. As Henri Theil
retired from the project soon after the European regional committee was
appointed, the main credit for what has been achieved should go to Jim
Durbin and Erich Lehmann.

%s5 ###
\section{Personal memories}
Over the years I spent many summers and two full semesters in Berkeley
and during these years I got to know Erich Lehmann pretty well. Erich
also visited the Netherlands a number of times and in 1985 Leiden
University awarded an honorary doctorate to him. This happens only
rarely and the last time a mathematical scientist became a doctor h.c.
was a century earlier in 1885 when the laureate was Thomas Stieltjes,
also a name well known to probabilists and statisticians. Erich has
given a humorous account of this event in the final chapter of
\citet{Leh08}.

A few months before the award ceremony a delicate problem landed on my
plate. When Leiden University was founded in 1575 after the siege of
Leiden by the Spanish was lifted, it adopted the motto Praesidium
Libertatis (Bulwark of Liberty). During World War II, the dean of the
law school lived up to this motto by protesting the firing of the
Jewish professors by the German occupation in a public address. Faculty
and students joined in this protest, the dean was arrested and the
Germans closed down the university. Given this history, I was not
surprised to receive a phone call from the rector who asked me whether
I was absolutely sure that Professor Lehmann had no Nazi taint in his
past. I told him that the Lehmann family left Germany immediately after
Hitler came to power in 1933 when Erich was 15, and that Erich spent a
year in Guam with the U.S. Air Force. The rector said this was all very
well, but he could have been a member of the Hitler Jugend while still
in Germany. I replied that anyone who knew Erich would find it very
difficult to picture him in a HJ uniform, but I also had to admit that
the rector had a legitimate reason for wanting to avoid any risk of a
possible scandal. Of course I~found the whole thing rather embarrassing
and I tried to ask in a roundabout way whether my friends in Berkeley
knew anything of relevance. Of course what happened was that nobody
seemed to know anything about this until somebody asked Erich's wife Julie
who provided a perfect answer. The Lehmann family was of Jewish origin. Erich's father had read
 Hitler's \textit{Mein Kampf} and, unlike many of his friends, took Hitler seriously. Moreover, as a prominent
lawyer, he had won a major lawsuit against some important Nazi. So when Hitler
came to power shortly after, he decided that the family should leave
 the country as quickly as possible, which they did. The
rector was convinced but I still felt embarrassed. However, I~need not
have worried. Much later Julie told me that both she and Erich had liked the fact that somewhere in the
world people still cared about these matters.

The day after Erich and Julie arrived for a week of
festivities surrounding the award ceremony, Erich had a fever that took
two days to control and allowed him to skip one of the nightly dinner
parties offered by the rector, the dean, the department, etc. After
returning to the U.S. he wrote to me: ``\textit{On the way back, I had one
insight into the proceedings of last week. Could the purpose of this
series of 10-course dinners be to dispose of the newly baked (doctor)
h.c.'s before they can do any damage to the reputation of the Acad.
Lugd. Batav. (i.e., Leiden University)?}''

In a similarly fanciful but somewhat more serious vein he continued with:
``\textit{Bill, neither of us tends to take things too seriously but
(\ldots) even at the risk of embarrassing you (I must tell you) of the
extraordinary, liberating effect (last week's ceremony) had on me.}

\textit{I grew up with great advantages and some serious disadvantages. The
most important of the latter was the relationship with my father, a
very outgoing, enormously successful lawyer who adored my equally
outgoing, delightful and original younger brother (\ldots) but found
little to admire in me, who was shy, reserved, and unadventurous. It is
not that he mistreated me---on the contrary, he was always very nice to
me---he just made it clear that I was a bit disappointing and that he
had no expectations of my ever amounting to anything. Over the years
the weight of his disapproval has lightened some, but it has remained
with me until last week when I had it out with him (the poor man has
been dead for 35 years). I challenged him: `Now admit that you were
wrong. In my quiet boring way I did make it after all---getting the
only mathematical Dr h.c. from Leiden in the 20th century is as big a
success as any man can wish for his son---so there!' And he admitted it}.''

In 1997 I spent two months in Berkeley. I was sorry to hear that Erich
and Julie would be in Princeton at that time, but happy when they
offered me the use of their apartment. For some reason they came back
earlier than planned, so after being there for a month, I prepared to
move out. However, they argued that since their apartment really
consisted of two apartments joined together by removing a wall, I could
still use the apartment on one side of the virtual wall for the
remainder of my stay. Of course I protested, but they refused to listen.

So for one month we shared their double apartment. Because Erich used
to get up at 5.30 a.m. to work on one of his books, he had finished the day's
work by the time I got up. We had a cup of coffee and a chat and I went
to the department. In the evening we had long discussions about
everything under the sun. In \citet{Leh08} Erich writes: ``\textit{In the
evenings, after dinner, we used to discuss the problems of the world
over a glass of genever, a supply of which he had brought from the
Netherlands}.'' We may not have solved all of these problems but it was
certainly fun.

In Erich's last letter to me in early 2009, there was another fanciful
dream. He wrote: ``\textit{Last night you made an appearance in my dream.
You introduced me to a colleague of yours who had a massive data set
that she needed to analyze and you asked me to assist her in this
enterprise. I reminded you that my facility with applied statistics
parallels that attributed to my mother's intelligence by one of her
teachers:}
\begin{center}
\begin{tabular}{ll}
\textit{Dumm geboren,} &\textit{Born stupid} \\
\textit{Nichts dazu gelernt} &\textit{Learned nothing since} \\
\textit{Und das noch vergessen} &\textit{And forgot even that.}
\end{tabular}
\end{center}

\textit{(A true story). But I said I'd do my best, which I hope you appreciate}.''

I wrote back that I was sorry to have disturbed his well-earned sleep
and that I liked the evaluation of his mother's intelligence, though
obviously false. Since neither Erich nor I are known for our applied
work, it was perhaps just as well that it was all a dream.

Erich Lehmann was a great man who touched many of our lives in a very
positive way. We are sad that he is no longer with us, but the good
memories remain.

%suskaldyti doi

% imsref loaded by smiklovaite, 2011-04-11 13:25:44
%
% imsref loaded by smiklovaite, 2011-04-11 15:22:19

\printaddresses


\begin{thebibliography}{17}
% BibTex style file: ims.bst, 2010-03-23
% Default style options (sort=0,type=number).
% Used options (sort=1,type=nameyear).

\bibitem[\protect\citeauthoryear{Albers, Bickel and van
  Zwet}{1976}]{AlbBicvan76}
\begin{barticle}[auto:STB|2011-03-03|12:04:44]
\bauthor{\bsnm{Albers},~\bfnm{W.}\binits{W.}},
  \bauthor{\bsnm{Bickel},~\bfnm{P.~J.}\binits{P.~J.}} \AND
  \bauthor{\bparticle{van} \bsnm{Zwet},~\bfnm{W.~R.}\binits{W.~R.}}
(\byear{1976}).
\btitle{Asymptotic expansions for the power of distribution free tests in the
  one-sample problem}.
\bjournal{Ann. Statist.}
\bvolume{4}
\bpages{108--156}.
\bid{mr={0391373}}
\end{barticle}
\endbibitem

\bibitem[\protect\citeauthoryear{Bickel and van Zwet}{1978}]{Bicvan78}
\begin{barticle}[auto:STB|2011-03-03|12:04:44]
\bauthor{\bsnm{Bickel},~\bfnm{P.~J.}\binits{P.~J.}} \AND
  \bauthor{\bparticle{van} \bsnm{Zwet},~\bfnm{W.~R.}\binits{W.~R.}}
(\byear{1978}).
\btitle{Asymptotic expansions for the power of distributionfree tests in the
  two-sample problem}.
\bjournal{Ann. Statist.}
\bvolume{6}
\bpages{937--1004}.
\bid{mr={0499567}}
\end{barticle}
\endbibitem

\bibitem[\protect\citeauthoryear{DeGroot}{1986}]{DeG86}
\begin{barticle}[mr]
\bauthor{\bsnm{DeGroot},~\bfnm{Morris~H.}\binits{M.~H.}}
(\byear{1986}).
\btitle{A conversation with {E}rich {L}. {L}ehmann}.
\bjournal{Statist. Sci.}
\bvolume{1}
\bpages{243--258}.
\bid{issn={0883-4237}, mr={0846003}}
\end{barticle}
\endbibitem

\bibitem[\protect\citeauthoryear{Feller}{1966}]{Fel66}
\begin{bbook}[mr]
\bauthor{\bsnm{Feller},~\bfnm{William}\binits{W.}}
(\byear{1966}).
\btitle{An Introduction to Probability Theory and Its Applications. {V}ol.
  {II}}.
\bpublisher{Wiley}, \baddress{New York}.
\bid{mr={0210154}}
\end{bbook}
\endbibitem

\bibitem[\protect\citeauthoryear{Fisher}{1925}]{Fis25}
\begin{barticle}[auto:STB|2011-03-03|12:04:44]
\bauthor{\bsnm{Fisher},~\bfnm{R.~A.}\binits{R.~A.}}
(\byear{1925}).
\btitle{Theory of statistical estimation}.
\bjournal{Proc. Cambridge Philos. Soc.}
\bvolume{22}
\bpages{700--725}.
\end{barticle}
\endbibitem

\bibitem[\protect\citeauthoryear{Fisher}{1956}]{Fis56}
\begin{bbook}[auto:STB|2011-03-03|12:04:44]
\bauthor{\bsnm{Fisher},~\bfnm{R.~A.}\binits{R.~A.}}
(\byear{1956}).
\btitle{Statistical Methods and Scientific Inference}.
\bpublisher{Oliver and Boyd}, \baddress{Edinburgh}.
\end{bbook}
\endbibitem

\bibitem[\protect\citeauthoryear{Hodges and Lehmann}{1970}]{HodLeh70}
\begin{barticle}[mr]
\bauthor{\bsnm{Hodges},~\bfnm{J.~L.}\binits{J.~L.} \bsuffix{Jr.}} \AND
  \bauthor{\bsnm{Lehmann},~\bfnm{E.~L.}\binits{E.~L.}}
(\byear{1970}).
\btitle{Deficiency}.
\bjournal{Ann. Math. Statist.}
\bvolume{41}
\bpages{783--801}.
\bid{issn={0003-4851}, mr={0272092}}
\end{barticle}
\endbibitem

\bibitem[\protect\citeauthoryear{Hoeffding}{1951}]{H1951}
\begin{bincollection}[auto:STB|2011-03-03|12:04:44]
\bauthor{\bsnm{Hoeffding},~\bfnm{W.}\binits{W.}}
 (\byear{1951}).
 \btitle{Optimum nonparametric tests}.
 In \bbooktitle{Proc. 2nd Berkeley Sympos. Math. Statist. Probab.}
 \bpages{83--92}.
 \bpublisher{Univ. California Press}, \baddress{Berkeley}.
\end{bincollection}
\endbibitem

\bibitem[\protect\citeauthoryear{Kallenberg et al.}{1984}]{Kal84}
\begin{bbook}[mr]
\bauthor{\bsnm{Kallenberg},~\bfnm{W.~C.~M.}\binits{W.~C.~M.}} \betal{et al.}
(\byear{1984}).
\btitle{Testing Statistical Hypotheses: Worked Solutions}.
\bseries{CWI Syllabi}
\bvolume{3}.
\bpublisher{Stichting Mathematisch Centrum, Centrum voor Wiskunde en
  Informatica}, \baddress{Amsterdam}.
\bid{mr={0778033}}
\end{bbook}
\endbibitem

\bibitem[\protect\citeauthoryear{Lehmann}{1953}]{Leh53}
\begin{barticle}[mr]
\bauthor{\bsnm{Lehmann},~\bfnm{E.~L.}\binits{E.~L.}}
(\byear{1953}).
\btitle{The power of rank tests}.
\bjournal{Ann. Math. Statist.}
\bvolume{24}
\bpages{23--43}.
\bid{issn={0003-4851}, mr={0054208}}
\end{barticle}
\endbibitem

\bibitem[\protect\citeauthoryear{Lehmann}{1959}]{Leh59}
\begin{bbook}[mr]
\bauthor{\bsnm{Lehmann},~\bfnm{E.~L.}\binits{E.~L.}}
(\byear{1959}).
\btitle{Testing Statistical Hypotheses}.
\bpublisher{Wiley}, \baddress{New York}.
\bid{mr={0107933}}
\end{bbook}
\endbibitem

\bibitem[\protect\citeauthoryear{Lehmann}{1997}]{Leh97}
\begin{barticle}[mr]
\bauthor{\bsnm{Lehmann},~\bfnm{E.~L.}\binits{E.~L.}}
(\byear{1997}).
\btitle{{T}esting statistical hypotheses: The story of a book}.
\bjournal{Statist. Sci.}
\bvolume{12}
\bpages{48--52}.
\bid{doi={10.1214/ss/1029963261}, issn={0883-4237}, mr={1466430}}
\end{barticle}
\endbibitem

\bibitem[\protect\citeauthoryear{Lehmann}{2008}]{Leh08}
\begin{bbook}[mr]
\bauthor{\bsnm{Lehmann},~\bfnm{E.~L.}\binits{E.~L.}}
(\byear{2008}).
\btitle{Reminiscences of a Statistician: The Company I Kept}.
\bpublisher{Springer}, \baddress{New York}.
\bid{doi={10.1007/978-0-387-71597-1}, mr={2367933}}
\end{bbook}
\endbibitem

\bibitem[\protect\citeauthoryear{Rao}{1961}]{Rao61}
\begin{bincollection}[mr]
\bauthor{\bsnm{Rao},~\bfnm{C.~Radhakrishna}\binits{C.~R.}}
(\byear{1961}).
\btitle{Asymptotic efficiency and limiting information}.
In \bbooktitle{Proc. 4th {B}erkeley {S}ympos. {M}ath. {S}tatist.
{P}robab.}
\bvolume{I}
\bpages{531--545}.
\bpublisher{Univ. California Press}, \baddress{Berkeley}.
\bid{mr={0133192}}
\end{bincollection}
\endbibitem

\bibitem[\protect\citeauthoryear{Reid}{1982}]{Rei82}
\begin{bbook}[mr]
\bauthor{\bsnm{Reid},~\bfnm{Constance}\binits{C.}}
(\byear{1982}).
\btitle{Neyman---From Life}.
\bpublisher{Springer}, \baddress{New York}.
\bid{mr={0680939}}
\end{bbook}
\endbibitem

\bibitem[\protect\citeauthoryear{van Dantzig}{1957}]{autokey3}
\begin{barticle}[auto:STB|2011-03-03|12:04:44]
\bauthor{\bparticle{van} \bsnm{Dantzig},~\bfnm{D.}\binits{D.}}
(\byear{1957}).
\btitle{Statistical priesthood II}.
\bjournal{Stat. Neerl.}
\bvolume{11}
\bpages{185--200}.
\end{barticle}
\endbibitem

\bibitem[\protect\citeauthoryear{Wilcoxon}{1945}]{Wil45}
\begin{barticle}[auto:STB|2011-03-03|12:04:44]
\bauthor{\bsnm{Wilcoxon},~\bfnm{F.}\binits{F.}}
(\byear{1945}).
\btitle{Individual comparisons by ranking methods}.
\bjournal{Biometrics Bulletin}
\bvolume{1}
\bpages{80--83}.
\end{barticle}
\endbibitem

\end{thebibliography}
\end{document}